\documentclass[12pt, amsfonts, epsfig]{amsart}
\usepackage{latexsym, amssymb, amsmath, epsfig, amscd, amsfonts}
 \newlength{\baseunit}               
 \newcount{\numlines}                
\newcommand{\Z}{\mathbb{Z}}

\newcommand{\Q}{\mathbb{Q}}
\newcommand{\R}{\mathbb{R}}
\newcommand{\C}{\mathbb{C}}

        \newfont{\hollow}{msbm10 scaled\magstep1}
        
        \newfont{\Bfmit}{eufm10 scaled\magstep1}
        \newcommand{\bfmit}[1]{\hbox{\Bfmit {#1}}}

\newcommand{\cL}{{\mathcal L}}

\newcommand{\cF}{{\mathcal F}}

\newcommand{\codim}{\operatorname{codim}}

\newcommand{\Pic}{\operatorname{Pic}}

\newcommand{\Cox}{\operatorname{Cox}}

\newcommand{\GL}{\operatorname{GL}}

\newcommand{\ev}{\operatorname{ev}}

\newcommand{\Hom}{\operatorname{Hom}}

\newcommand{\Spec}{\operatorname{Spec}}

\newtheorem{thm}{Theorem}[section]
\newtheorem{question}[thm]{Question}

\newtheorem{cor}[thm]{Corollary}

\newtheorem{prop}[thm]{Proposition}

\newtheorem{conj}[thm]{Conjecture}

\theoremstyle{definition}

\newtheorem{prob}[thm]{Problem}

\theoremstyle{remark}

\newcommand{\lremind}[1]{{}}
\newcommand{\bremind}[1]{{}}

\newcommand{\cut}[1]{}

\begin{document}
\pagestyle{plain} \title{{ \large{Geometric Invariant Theory and
Birational Geometry} } }
\author{Yi Hu}

\address{Department of Mathematics, University of Arizona, Tucson,
AZ 86721, USA}

\email{yhu@math.arizona.edu}

\address{Center for Combinatorics, LPMC, Nankai University, Tianjin 300071, China}


\maketitle

{\parskip=12pt 

\begin{abstract}
In this paper we will survey some recent developments in the last
decade or so on variation of Geometric Invariant Theory and its
applications to Birational Geometry such as the weighted Weak
Factorization Theorems of nonsingular projective varieties and
more generally projective varieties with finite quotient
singularities. Along the way, we will also mention some progresses
on birational geometry of hyperK\"ahler manifolds as well as
certain open problems and conjectures.
\end{abstract}


\section{Introduction}

Quotients of projective varieties by reductive algebraic groups
arises naturally in many situations. The existence of many
important moduli spaces, for example, is proved by expressing them
as quotients. There are several quotient theories, among them, the
Mumford geometric invariant theory (GIT) is the systematic one. A
basic observation there is that many moduli functors can be, at
least coarsely, represented by quotient varieties in the sense of
GIT.

 A remarkable discovery in the last decade is the deep
connection and fruitful interaction between GIT and birational
geometry (e.g.,
\cite{DH},\cite{Hu2004},\cite{HuKeel2},\cite{HuKeel},\cite{Th},
\cite{Wlodarczyk00}).

Moduli spaces, or GIT quotients are not unique in general. Hence a
natural fundamental problem in this theory is to parameterize all
the GIT quotients and describe exactly how a quotient changes as
one varies the underlying parameter.

This fundamental question was first raised and answered
 in the joint work \cite{DH}, which overlaps with the paper \cite{Th}.
Here, one may wonder why the founder of the theory, in his famous
work \cite{GIT},  did not mention such a fundamental question.
But, it is really not that surprising given that the higher
dimensional birational geometry witnessed a great advance only
when Mori revolutionized the field
 many years after GIT was invented.
Our solution to this fundamental problem is a successful mixture
of the  great works of the two Fields medalists.

Roughly speaking, the solution given by the theory of Variation of
Geometric Invariant Theory quotients (VGIT) can be summarized as
follows. For a given action, there are finitely many GIT quotients
and they are parameterized by some natural chambers in the
effective equivariant  ample cone. Moreover, when crossing a wall
in the chamber system, the corresponding quotient undergoes a
birational transformation similar to a Mori flip.

This, for the first time, provided a significant link between VGIT
and Mori program. Later (\cite{HuKeel2},
\cite{HuKeel},\cite{Hu2004}), we found that the connection between
the two is in fact deeper. To be short, we discovered that not
only VGIT wall-crossing transformations provide special examples
of transformations in birational geometry, but remarkably all the
birational transformations in birational geometry can also be
decomposed as sequences of  VGIT wall-crossing transformations.
This establishes a useful philosophy that began  in the paper
\cite{HuKeel2} and reconfirmed in \cite{HuKeel} and \cite{Hu2004}:
Birational geometry is a special case of Variation of Geometric
Invariant Theory.

The philosophy has two versions: global and local.

The global version says that for a given $\Q$-factorial projective
variety and all its Mori flip images, they all can be realized as
GIT quotients of certain action. If this holds, we will call the
projective variety a {\it Mori dream space}. Not all
$\Q$-factorial projective varieties are Mori dream space. In
\cite{HuKeel2}, we give a criterion for a Mori dream space in
terms of the Cox section ring.

The local version says that any two smooth projective varieties,
or more generally, any two projective varieties with finite
quotient singularities related by a birational morphism can be
realized as GIT quotients of a smooth projective variety by a
reductive algebraic group. This turns out to be always true. For
the smooth case, the GIT realization of a birational morphism was
stated and proved  by Hu and Keel in \cite{HuKeel}. We were
informed that it may also follow from Wlodarczyk's paper
\cite{Wlodarczyk00} (see also \cite{AKMW}). A consequence of the
above is the GIT Weak Factorization Theorem which asserts that any
two birational smooth projective varieties are related by a
sequence of VGIT wall-crossing birational transformations. Note
here that the (stronger) Weak Factorization Theorems for the
smooth cases were proved by Wlodarczyk et. al., see \cite{AKMW}
and the reference therein. A proof for singular cases was given in
\cite{Hu2004}. The singular version of the Weak Factorization
Theorems asserts that any two birational projective varieties with
finite quotient singularities are related by a sequence of VGIT
wall-crossing birational transformations (\cite{Hu2004}). Here, by
 a VGIT wall-crossing birational transformation, we mean the
 birational transformation as described in Theorem 4.2.7 of \cite{DH}.

Besides the GIT approach, another way of compactification in
moduli problem is Chow quotient. This approach to moduli problems
is different and may be harder than GIT, but it has some
advantages and potentially has several significant applications.
Understanding the boundary of Chow quotient is the key to this
approach. In \cite{Hu2003b}, we give a characterization regarding
the boundaries of Chow quotient, using a very intuitive and
computable method, namely the
Perturbation-Translation-Specialization relation. This
characterization, for example, can be used to give symplectic and
topological interpretations of Chow quotient.

Turning away from the birational geometry of general varieties, we
now would like to focus on special varieties. HyperK\"ahler
Varieties form a sepcial and very important class of Calabi-Yau
manifolds. Its birational geometry has received some considerable
attentions recently (\cite{BHL}, \cite{miyaoka}, \cite{HY},
\cite{Fu}, \cite{WW}, among others). Being very restrictive in
structures, birtaional transformations among HyperK\"ahler
Varieties are among simplest kinds. It was proved in \cite{BHL},
coupled with an improvement of \cite{WW}, that any two birational
HyperK\"ahler varieties of dimension 4 are related by a sequence
of Mukai elementary transformations. In higher dimensions, Hu and
Yau (\cite{HY}) obtained several structure theorems for symplectic
contractions, completing important initial steps toward the Hu-Yau
conjecture on decomposing symplectic birational maps. We remark
here that some of our results in this direction overlap with the
work of \cite{miyaoka}

\section{Moduli, GIT, and Chow Quotient}

\subsection{What are Moduli Spaces?}
A moduli space is, roughly speaking, a parameter space for
equivalence classes of certain geometric objects of fixed
topological type. Depending on purposes, there are two types of
moduli: fine moduli and coarse moduli. The former often does not
exist as a (quasi-) projective scheme but lives naturally as a
stack, while the latter frequently exists as a (quasi-) projective
scheme.

A naive way to think of a moduli space $M$ is
$$M =\{\rm Geometric \; Objects \} / \sim$$
as a collection of the geometric objects of the interest modulo
equivalence relation. This way of thinking is "coarse". Such a
moduli space $M$ often exists as a projective scheme.

To describe fine moduli, one needs moduli functor, a categorical
language for the problem. A moduli functor is covariant functor
$$\cF: \{\rm schemes\} \rightarrow \{\rm sets\}$$
which sends any scheme $X$ to a family of the geometric objects
parameterized by $X$. The moduli space $M$ is fine, or the moduli
functor $\cF$ is represented by $M$ if there is a universal family
over $M$
$$U \rightarrow M$$ such that for any scheme $X$ and a family of
the geometric objects
$$Z \rightarrow X$$
parameterized by $X$, it corresponds to a morphism $f: X
\rightarrow M$ so that the family $Z \rightarrow X$ is the
pullback of the universal family via the morphism $f$
\begin{equation*}
\begin{CD}
Z=f^*U @>>> U\\
@VVV   @VVV\\
X@>{f}>> M
\end{CD}
\end{equation*}

Fine moduli in general does not exist if some objects possess
nontrivial automorphisms. Therefore we are somehow forced to
consider coarse moduli if we prefer to work on projective schemes.

Now let us return to a coarse moduli as a parameter space of
geometric objects
$$M =\{\rm Geometric \; Objects \} / \sim.$$
This way of thinking, as it stands,  naturally lead to a quotient
theory. First we denote the collection of the objects as
$$X = \{\rm Geometric \; Objects \}$$ and then we would introduce a natural group action
$$G \times X \rightarrow X$$ such that
two geometric objects $x, y$ are equivalent if and only if they,
as points in $X$, are in the same group orbit:
$$x \sim y \;{\rm if \; and \; only \; if}\; G \cdot x = G \cdot
y.$$

This is what the general idea is. However, taking quotients in
algebraic geometry can be very subtle. Indeed, there are several
theories about this:
\begin{itemize}
\item The Hilbert-Mumford Geometric Invariant Theory;
\item Chow
quotient varieties  \item Hilbert quotient varieties; \item
Artin, Koll\'ar, Keel-Mori quotient spaces.
\end{itemize}

We will only focus on the first two: GIT quotients and Chow
quotients.

\subsection{What are GIT quotients?}

\subsubsection{A toy example}

To give the reader some intuitive  ideas about GIT quotients, let
us {\it informally} consider a simple, yet quite informative
``toy'' example.

Let $G = \C^*$ act on ${\mathbb P}^2$ by $$ \lambda \cdot [x:y:z]
= [\lambda x: \lambda^{-1} y : z].$$ Consider a map $\Phi:
{\mathbb P}^2 \rightarrow \R$ given by $$\Phi([x:y:z]) =
\frac{|x|^2 - |y|^2}{|x|^2 + |y|^2 + |z|^2}.$$ This is the
so-called moment map for the induced symplectic $S^1$-action with
respect to the Fubibi-Study metric. Its image is the interval
$[-1,1]$.

\bigskip

\begin{picture}(20, 9)
\put(4,.2){ \psfig{figure=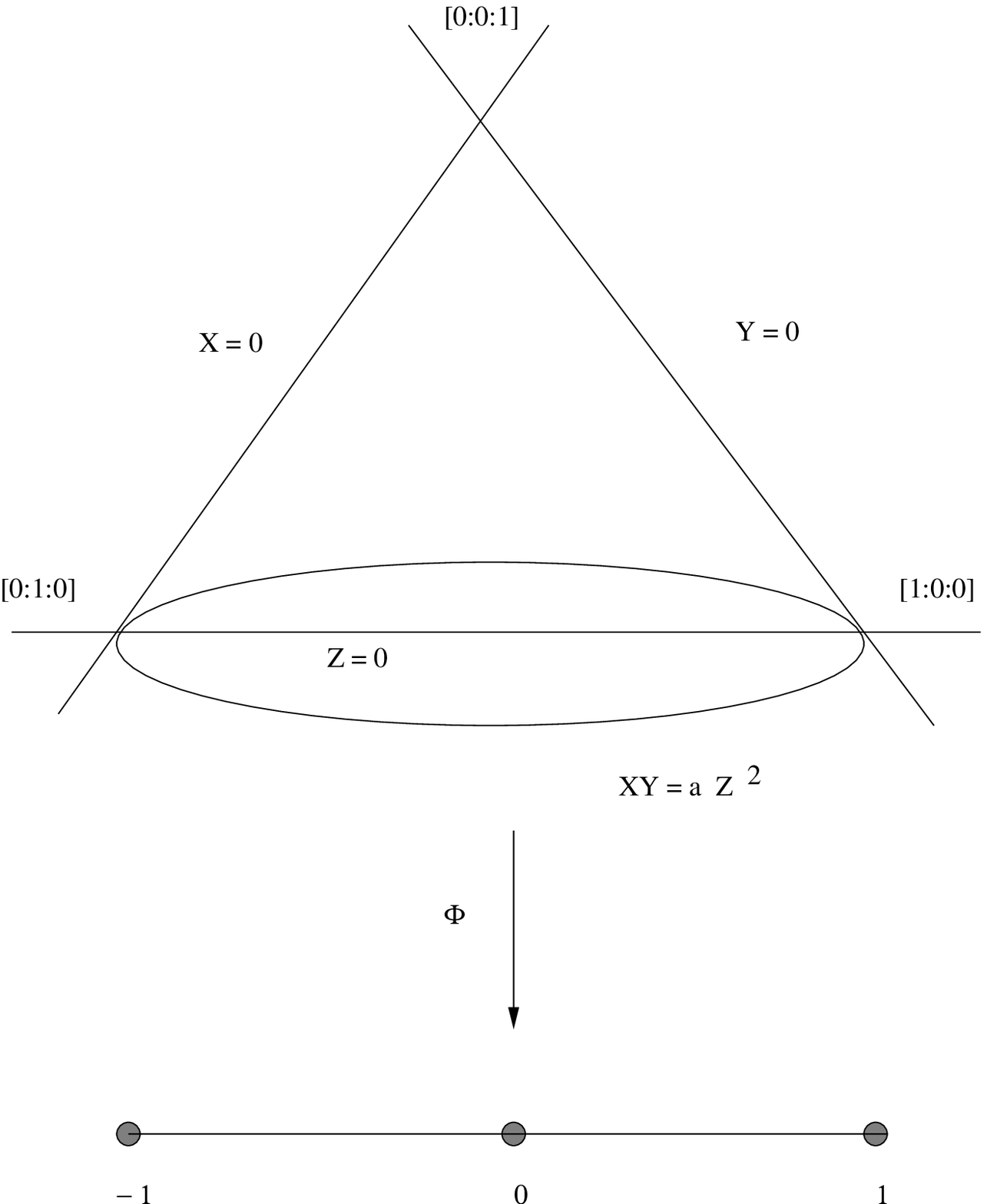, height = 8.5cm,width = 8cm}
}
\end{picture}

\bigskip
\centerline{Figure 1. \ Conics}





The $\C^*$ orbits are classified as follows.  (See Figure 1 for an
illustration.)
\begin{itemize}
\item Generic $\C^*$-orbits are conics $XY = a Z^2$ minus two
points [1:0:0] and [0:1:0] for $a \ne 0, \infty$.  We denote these
orbits by $(XY = a Z^2)$.

 The moment map image of the orbit $(XY =
a Z^2)$ is $(-1,1)$.

 \item Other 1-dimensional orbits
are the three coordinate lines $X=0$, $Y=0$, and $Z=0$ minus the
coordinate points on them. We denote these orbits by $(X=0)$,
$(Y=0)$, and $(Z=0)$.

The moment map images of the orbits $(X=0)$, $(Y=0)$, and $(Z=0)$
are $(-1,0), (0,1)$, and $(-1,1)$, respectively.

\item Finally, the fixed points are the three coordinate points,
[1:0:0], [0:1:0], and [0:0:1].

Their moment map images are $1, -1,$ and $0$, respectively.
\end{itemize}

The reader can easily check that the ordinary topological orbit
space ${\mathbb P}^2/\C^*$ is a nasty non-Hausdorff space. For
example, the orbits $(X=0)$ and $(Y=0)$ are in the limits of the
generic orbits $(XY = a Z^2)$ when $a$ approaches 0, thus can not
be separated.

The idea of GIT is to find some open subset $U \subset {\mathbb
P}^2$, called the set of semistable points, such that a good
quotient in a suitable sense, $U /\!/ \C^*$, exists. In this
particular example, such open subsets are selected as follows.

The moment map $\Phi$ has three critical values $-1, 0,$ and $1$
which divide the interval into two top chambers $[-1,0]$ and
$[0,1]$, and three 0-dimensional chambers $\{-1\}, \{0\}, \{1\}$.
Each chamber $C$ defines a GIT stability:  a point $[x:y:z]$ is
semi-stable with respect to the chamber $C$ if $$C \subset
\Phi(\overline{\C^* \cdot [x:y:z]}),$$ and it is stable if the
(relative) interior $C^\circ $ of $C$ satisfies
$$C^\circ \subset
\Phi(\C^* \cdot [x:y:z])\; \hbox{and}\; \dim \C^* \cdot [x:y:z] =
1.$$ (For a reference for this, see for example, \cite{HuMIT}.)
Thus, for example, the orbit $(X=0)$ is stable with respect to
$[-1,0]$, unstable with respect to $[0,1]$; while the orbit
$(Y=0)$ is stable with respect to $[1,0]$, unstable with respect
to $[-1,0]$. But, $(X=0)$, $(Y=0)$, and $[0:0:1]$ are all
semi-stable with respect to the chamber $\{0\}$. Finally, observe
that the generic orbits, namely the conics $(XY = a Z^2)$ are
stable with respect to any chamber.

A general philosophy of GIT quotient is that it should
parameterize orbits that are {\it closed} in the semi-stable
locus. Here in this example, the GIT quotient $X_{[-1,0]}$ defined
by the chamber $[-1,0]$ parameterizes $(XY = a Z^2)$ $(a \ne 0,
\infty)$, $(Z=0)$, and $(X=0)$. Thus $X_{[-1,0]}$ is isomorphic to
${\mathbb P}^1$ with $(Z=0)$ and $(X=0)$ serve as $\infty$ and 0,
respectively. Likewise, the GIT quotient $X_{[0,1]}$ defined by
the chamber $[0,1]$ parameterizes $(XY = a Z^2)$, $(Z=0)$, and
$(Y=0)$. And, the GIT quotient $X_{\{0\}}$ defined by the chamber
$\{0\}$ parameterizes $(XY = a Z^2)$, $(Z=0)$, and $[0:0:1]$.

All these quotients are isomorphic to ${\mathbb P}^1$, and hence
they are all isomorphic to each other. This is very special
because the dimension is too low (namely 1) to allow any
variation. In general, they should be quite different and are only
birational to each other (see \S 2.1). This was a very decisive
observation, and the determination to investigate the most general
relation among different GIT quotients led to some significant
discoveries and applications.

\subsubsection{GIT quotients in general} In general, GIT quotients
are constructed as follows. Consider an action $$G \times X
\rightarrow X$$ with $X$ a projective variety and $G$ a reductive
algebraic group. To define a GIT quotient, we need to make some
choices
\begin{equation*}
\begin{CD}
G \times L @>>> L \\
@VVV   @VVV\\
G \times X @>>> X
\end{CD}
\end{equation*}
where $L$ is an ample line bundle over $X$ and the action on $X$
is lifted to $L$ such that the induced map
$$g: L_x \rightarrow L_{g \cdot x}$$
is a linear map for any $x \in X$ and $g \in G$. Here $L_x$
denotes the fiber of $L$ over the point $x$. Such a device is
called a linearized ample line bundle or simply a linearization,
which we still denote by $L$ for simplicity. Because of the lifted
action, the space of sections $\Gamma(X,L)$ becomes a
representation of $G$. That is, the linearization induces a
$G$-action on $\Gamma (X, L)$ as follows:
for any $s \in \Gamma
(X, L)$ and $g \in G$
$$ ^g s (x) = g \cdot s (g^{-1} \cdot x).$$
We say  a section $s$ is $G$-invariant if $ ^g s = s$ for all $g
\in G$. This simply means that
$$s(g \cdot x) = g \cdot s(x)$$
for all $g \in G$ and $x \in X$. We will denote the fixed point
set of $\Gamma(X,L)$ by $\Gamma^G(X,L)$.

The linearization $L$ singles out a distinguished Zariski open
subset $X^{ss}(L)$ of $X$ as follows:
$$X^{ss}(L) =\{x \in X | \; {\rm there \; is \;} s \in
\Gamma^G(X,L) \;{\rm such \; that\;} s(x) \ne 0 \}.$$ Points in
$X^{ss}(L)$ are called semistable points with respect to the
linearization $L$. $X^{ss}(L)$ contains a generally smaller subset
$$X^s = \{ x \in X^{ss}(L) |\;G \cdot x \;
{\rm is \; closed \; in\;} X^{ss}(L) \; {\rm and\;} \dim G \cdot x
= \dim G \}$$ whose points are called stable points w.r.t the
linearization $L$.

The quotient of GIT is different than topological quotient in that
it is not orbit space in general. GIT introduces the following
{\it new} relation in $X^{ss}(L)$: $$x \sim y \;\hbox{if and only
if }\; \overline{G \cdot x} \cap  \overline{G \cdot y} \cap
X^{ss}(L) \ne \emptyset.$$ Topologically, the GIT quotient is
$$X^{ss}(L)/\!/G = X^{ss}(L)/\sim.$$ From the definition, one
checks that $X^{ss}(L)/\!/G$ parameterizes the closed orbits in
$X^{ss}(L)$, and the equivalence relation on $X^s$ is the same as
the orbit relation, hence the orbit space $X^s/G$ is contained in
$X^{ss}(L)/\sim.$

A basic theorem of the geometric invariant theorem says

\begin{thm}{\rm (\cite{GIT})}
The GIT quotient $X^{ss}(L)/\!/G$ exists as a projective variety
and contains the orbit space $X^s/G$ as a Zariski  open subset.
\end{thm}

\subsection{Relation with Symplectic Reduction}

There is a beautiful link between Geometric Invariant Theory and
Symplectic Geometry.

Consider again an action $G \times X \rightarrow X$ of an
algebraic reductive group on a smooth projective variety $X$, and
let $L$ be a linearized ample line bundle. The line bundle induces
 an integral symplectic form $\omega$ on $X$, the linearization
amounts to a Hamiltonian symplectic action of a compact form $K$
of $G$ on the symplectic manifold $(X, \omega)$ and an equivariant
differentiable map, called moment map
$$\Phi: X \rightarrow {\bfmit k}^*$$
where ${\bfmit k}$ is the Lie algebra of  the compact Lie group
$K$. The moment map is a solution of the following differential
equation:
$$d (\Phi \cdot a) = i_{\xi_a} \omega$$
where $a \in {\bfmit k}$, $\xi_a$ is the vector field generated by
$a$, and $\Phi \cdot a$ is the component of $\Phi$ in the
direction of $a$.

The symplectic reduction or symplectic quotient by definition is
$$\Phi^{-1}(0)/K.$$ It inherits a symplectic form from $\omega$ on
its smooth part. A beautiful theorem, due to contribution from
many people, is that $\Phi^{-1}(0)$ is included in $X^{ss}(L)$ and
this inclusion induces a homeomorphism
$$\Phi^{-1}(0)/K \cong X^{ss}(L)/\!/G.$$

This relation can be extended to include K\"ahler quotients
(\cite{DH}).

\subsection{Relative GIT and Universal Moduli Spaces}

Stable locus behave nicely under an equivariant morphism.

Let  $f: Y \rightarrow X$ be a projective morphism between two
smooth projective varieties $X$ and $Y$, acted upon by two
reductive algebraic group $G'$ and $G$, respectively. Assume
further that we have an epimorphism $\rho: G' \rightarrow G$ with
respect to which $f: Y \rightarrow X$ is equivariant. Let $G_0$ be
the kernel of $\rho$. Then we have

\begin{thm} {\rm (--,\cite{Hu96})} Let $L$ and $M$ be two
linearized ample line bundle over $X$ and $Y$, respectively. Then
for sufficiently large $n$, we have
\begin{enumerate}
\item $Y^{ss}(f^*L^n \otimes M) \subset Y_{G_0}^{ss}(M) \cap
f^{-1}(X^{ss}(L))$; \item $Y^s(f^*L^n \otimes M) \supset
Y_{G_0}^s(M) \cap f^{-1}(X^s(L))$.

Assume in addition that $X^{ss}(L) = X^s (L)$, then \item
$Y^{ss}(f^*L^n \otimes M) = Y_{G_0}^{ss}(M) \cap
f^{-1}(X^{ss}(L))$; \item $Y^s(f^*L^n \otimes M) = Y_{G_0}^s(M)
\cap f^{-1}(X^s(L))$.
\end{enumerate}
\end{thm}

When $G' = G$, this is the main theorem of Reichstein
(\cite{Reichstein89}).

The theorem has best applications when the semistable locus and
stable locus coincide on the base (i.e., the cases (3) and (4)).
For example, this can be easily applied to the universal moduli
space of semistable coherent sheaves over the moduli space
$\overline{M_g}$ of stable curves of genus $g \ge 2$, which was
constructed  by Pandharipande (\cite{Pandharipande94a}), using
relative GIT. This is re-done in \cite{Hu96}, combining the above
theorem and the work of Simpson (\cite{Simpson94}).

\bigskip
\subsection{Chow quotients}

In the ``toy" example, observe that the lines $X=0$ and $Z=0$,
which are of degree 1,  have different homology classes than the
conic orbits $XY = a Z^2$, which are of degree 2. But the GIT
quotient $X_{[-1,0]}$ parameterizes these orbits of different
homology classes.  Even worst, the two orbits $(X=0)$ and $(Y=0)$
 are all identified with the closed orbit $[0:0:1]$ in the quotient $X_{\{0\}}$,
  and the orbit $[0:0:1]$ even has smaller dimension than the dimension of
the generic conic orbits. These are not desirable or suitable for
moduli problems.

This kind of problem can however be overcome by considering Chow
quotient which takes completely different approach.

Return to our ``toy example", to obtain the Chow quotient, we
first consider the closures of the generic $\C^*$-orbits, $XY = a
Z^2$ ($a \ne 0, \infty$), and then looks at all their possible
degenerations. When $a=0$, we get the degenerated conic $XY=0$,
two crossing lines; and when $a = \infty$, we obtain $Z^2 =0$, a
double line. They all have the same homology classes (degree 2).
And the Chow quotient is {\it the space of all $\C^*$-invariant
conics}. (See Figure 1.) Each point of the Chow quotient
represents an invariant algebraic cycle. In this case, the generic
cycles are $[XY = a Z^2]$ ($a \ne 0, \infty$), and the special
cycles are $[X=0] + [Y=0]$ ($a=0$), and $2 [Z=0]$ ($a=\infty$).

From  the above example, we see that Chow quotient parameterizes
cycles of generic orbit closures and their toric degenerations
which are {\it certain sums of orbit closures of dimension $\dim
G$}. We will call these cycles {\it Chow cycles} or {\it Chow
fibers}.

So, when do two arbitrary points belong to the same Chow cycle?

Consider the example again. We have that $[0:y:z]$ and $[x:0:z]$
($xyz \ne 0$) belong to the same Chow cycle $XY=0$. To get
$[x:0:z]$ from $[0:y:z]$, we first {\it perturb} $[0:y:z]$ to a
general position $$\varphi (t) = [tx:y:z] (t \ne 0), $$ then {\it
translate} it by $g(t) = t^{-1} \in \C^* $ to
$$g (t) \cdot \varphi
(t) = [x:ty:z],$$ and then $g (t) \cdot \varphi (t)$ {\it
specializes} to $[x:0:z]$ when $t$ goes to 0. We will call this
process {\it perturbing-translating-specializing} (PTS). It turns
out this simple relation holds true in general. That is, we prove
in general that two points $x$ and $y$ of $X$, with
$$\dim G \cdot
x = \dim G \cdot y = \dim G,$$ belong to the same Chow cycle if
and only if $x$ can be perturbed (to general positions),
translated along $G$-orbits (to positions close to $y$), and then
specialized to the point $y$.¡¡

We will return to this in \S 4.

\section{Variation of GIT and Factorization Theorem}

\subsection{Variation of GIT Quotients}
From the earlier discussion, we see that  GIT quotients are not
unique. Hence a natural fundamental problem in this theory is to
parameterize all the GIT quotients and describe exactly how a
quotient changes as one varies the underlying parameter. Here we
will explain  the solution to this problem as given  in \cite{DH}.
To motivate it, we begin with an example.

\subsubsection{An Example}

In the ``toy" example, we worked out all the GIT and Chow
quotients, but they are all isomorphic to ${\mathbb P}^1$, the
unique compactification of $\C^*$,  because we insist an example
that are very simple to describe. Here it should be fair to at
least point out to the reader a workable example where a
nontrivial wall crossing phenomenon and different quotients do
occur. Hence we take the liberty to include the following example
with  details left to the reader.

Consider the action of $\C^*$ on ${\mathbb P}^3$ by
$$ \lambda \cdot [x:y:z:w] = [\lambda x: \lambda y : \lambda^{-1}z: w].$$
The moment map is
$$\Phi([x:y:z:w]) = \frac{|x|^2 + |y|^2 - |z|^2}{|x|^2 + |y|^2 + |z|^2 + |w|^2}.$$

The image $\Phi(X)$ is $[-1,1]$ with three critical values $-1, 0,
1$.  So, we consider the level sets $\Phi^{-1}(-\frac{1}{2}),
\Phi^{-1}(0), \Phi^{-1}(\frac{1}{2})$. In Figure 2, we illustrate
the {\it real parts}
$$\Phi^{-1}(-\frac{1}{2})_\R, \Phi^{-1}(0)_\R,
\Phi^{-1}(\frac{1}{2})_\R$$ of the level sets restricted to $\C^3
\subset {\mathbb P}^3$ (the $\C^3$ is defined by setting $w= 1$).
It turns out this real picture preserves all the topological
information we need.

To understand this picture, note that the real points of $S^1$ are
$$S^1 \cap \R =  \{-1, 1\} = \Z_2.$$ So, the real parts of the
symplectic quotients (or the GIT quotients)
$$\Phi^{-1}(-\frac{1}{2})/S^1, \Phi^{-1}(0)/S^1,
\Phi^{-1}(\frac{1}{2})/S^1$$  are
$$\Phi^{-1}(-\frac{1}{2})_\R/\Z_2, \Phi^{-1}(0)_\R/\Z_2,
\Phi^{-1}(\frac{1}{2})_\R/\Z_2.$$

\begin{picture}(15, 9)
\put(4,1){ \psfig{figure=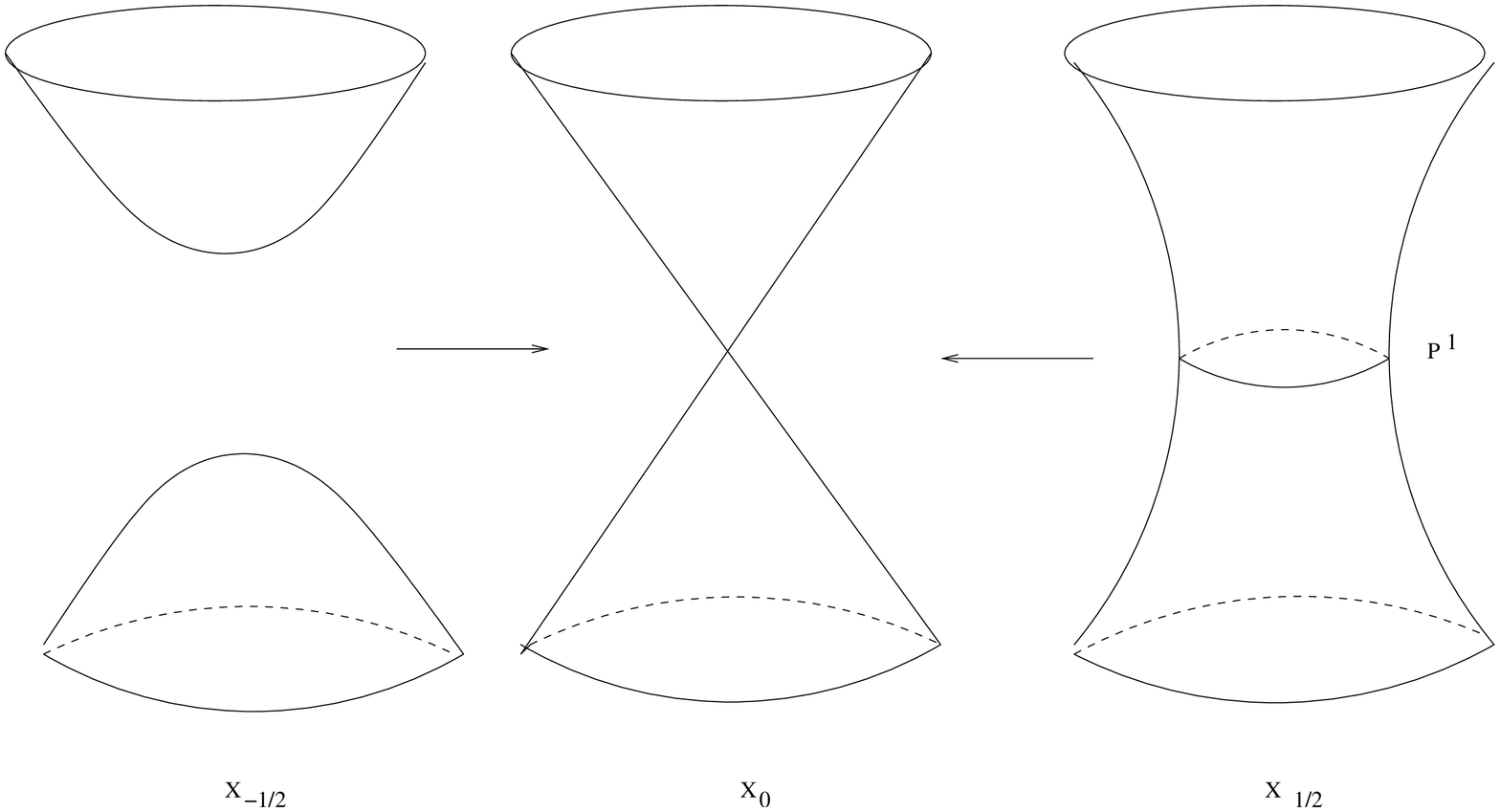,height = 7.5cm,width = 8cm} }
\end{picture}

\centerline{Figure 2.  Wall-Crossing Maps}

\bigskip

Note that $\Z_2$ acts on each level set by identifying the lower
part with the upper part. Hence, the quotient can be naturally
identified with the upper part.

There are two natural collapsing maps,  as shown in the picture.

 Now
observe that the left map happens to be an isomorphism, but the
right map is a (real) blowup along the origin so that the special
fiber is $\R {\mathbb P}^1$. Now complexifying this picture and
compactifying the results, we obtain three quotients: $X_{[-1,0]}
\cong {\mathbb P}^2$, $X_{\{0\}} \cong {\mathbb P}^2$, and
$X_{[0,1]}$ isomorphic to the blowup of ${\mathbb P}^2$ along a
point.

The reader may try to classify all the generic $\C^*$-orbits and
study their degenerations. He can verify that the Chow quotient is
also isomorphic to the blowup of ${\mathbb P}^2$ along a point.

This example suggests that when crossing a critical value of a
moment map, the corresponding GIT quotient changes by a blowdown
followed by a blowup. As it turns out, this simple wall-crossing
phenomenon exhibited in this simple example reveals a general
nature of GIT quotients and birational geometry.

\subsubsection{The Main Theorem of VGIT \cite{DH}}

 Some of the main results
of {\bf  \cite{DH}} may be roughly summarized as

\begin{thm} {\rm (Dolgachev-Hu, {\bf \cite{DH}})}
\label{vgit} Let $X$ be a smooth projective variety acted upon by
an algebraic reductive group $G$. Then there is a convex cone
$C^G(X)$, the $G$-ample cone,  which admits natural finite wall
and chamber structure such that {\it under very general
conditions}, we have
\begin{enumerate}
\item all GIT quotients are parameterized by rational points of
the cone $C^G(X)$. (Other points of $C^G(X)$ parameterize K\"ahler
quotients;) \item points in the same chamber give rise to
identical quotient; \item when crossing a wall, the quotient
undergoes a birational transformation similar to a Mori flip, VGIT
wall-crossing flip.
\end{enumerate}
\end{thm}

The main part of this theorem was also proved by Thaddeus
(\cite{Th}). The special case of this theorem when the group $G$
is Abelian was done
 earlier in \cite{BP} and independently in the author's MIT thesis (\cite{HuMITthesis}).
Theorem \ref{vgit} can be applied to a number of moduli problems,
some of which, e.g., moduli of parabolic bundles over a curve,
were done in the author's joint works (\cite{BH}). Many other nice
applications, especially to Donaldson invariants, appeared in the
works of several other mathematicians. The entire VGIT project,
the Abelian and the general cases, followed an observation of
Goresky and MacPherson, made in 80's (\cite{GM}), which alluded
that GIT morphisms should have implications for birational
geometry.

\subsection{VGIT and the Mori Program}

Theorem \ref{vgit} shows that VGIT provides natural examples of
factorizing general birational transformations into sequences of
simple ones. In my joint works \cite{HuKeel2} and \cite{HuKeel},
we found that the connection between the two is in fact deeper.
These two works were motivated, apart from others, by a natural
question that basically asks the converse of Theorem \ref{vgit}:

\begin{question}
Given any {\it general} birational map, can it be realized as a
sequence of VGIT wall-crossing flips?
\end{question}

In the special case when the birational map is made of a single
flip, this was answered by M. Reid, see the "tabernacle" lecture
by Reid, see also Theorem 1.7 of \cite{Th} (these were pointed out
by Abramovich). Note here that the existence of flips is not known
in higher dimensions.

In general, the answer may be divided into two cases: (1). VGIT
$\Leftrightarrow$ Mori: the global version; (2) VGIT
$\Leftrightarrow$ Mori: the local version. One direction of VGIT
$\Leftrightarrow$ Mori is Theorem \ref{vgit}. We only need to
explain the other.

We begin with the global version. It roughly goes as follows.
 Consider a ${\mathbb Q}$-factorial projective variety $X_0$
and all of its Mori flip images $\{X_i\}_i$.
The best thing one can possibly hope here is
that
there exists a group action $G \times \widetilde{X} \rightarrow \widetilde{X}$
such that all $\{X_i\}_i$ are GIT quotients of the action and all the flips
are VGIT wall-crossing flips.
If this is indeed the case, we will call $X_0$ a {\it  Mori dream space}.
We have a criterion for a Mori dream space.

\begin{thm} {\rm (Hu-Keel, {\bf \cite{HuKeel2}})} Let $X$ be a ${\mathbb Q}$-factorial
quasi-projective variety with $\Pic (X)_{\mathbb Q} = N^1 (X)$.
Then $X$ is a Mori dream space if and only if the Cox ring
$$\Cox (X) := \bigoplus_{(n_1, \ldots, n_r) \in {\mathbb N}^r}
 H^0(X, L_1^{n_1} \otimes \cdots \otimes L_r^{n_r})$$ is
finitely generated, where $L_1, \ldots, L_r$ is a basis of $N^1 (X)$.
In this case,
$X$ is naturally a quotient of the affine variety $\Spec (\Cox (X))$
by the group $H= \Hom({\mathbb N}^r, {\mathbb C}^*)$.
\end{thm}

Toric varieties are special cases of Mori dream spaces when $\Cox
(X)$ are polynomial rings (\cite{Cox}).  Other examples of Mori
dream spaces include quotients of Grassmanians by maximal tori, or
equivalently the moduli spaces of configuration of points in
projective spaces, and GIT quotients of some affine varieties. Of
course, most varieties cannot be Mori dream spaces, for example,
when the ample cone of $X$ is not polyhedral.

\subsection{VGIT and Weak Factorization Theorem}

However, the local version,  as we will explain now,
 is basically always valid. This version basically asks if  VGIT $\Leftrightarrow$ Mori
is true for some part of birational transformations. Indeed,
instead of putting all of birational transformations of $X_0$ in
the framework of just one group action (the Mori dream case), we
may use possibly different group actions for different
 birational morphisms.

 To this end, we have

\begin{thm} ({\rm Hu-Keel}, {\bf \cite{HuKeel}}. {\rm Wlodarczyk}, {\bf \cite{Wlodarczyk00}})
\label{fact} Let $X \rightarrow X'$ be a birational morphism
between two smooth projective varieties.
Then there is a ${\mathbb C}^*$- smooth projective variety $W$
such that $X$ and $X'$ are geometric GIT quotients of two open
subsets of $W$ by ${\mathbb C}^*$. Consequently, by VGIT,
 any birational map $f: X ---> Y$ between two smooth projective varieties
 can be factorized as a sequence of GIT wall-crossing flips.
\end{thm}

Here, by a GIT wall-crossing flip, we mean the birational
transformation as described in Theorem 4.2.7 of \cite{DH}. The GIT
realization of the birational morphism $X \rightarrow X'$ was
stated and proved in \cite{HuKeel}. We were informed that the same
may also be done with Wlodarczyk's construction in
\cite{Wlodarczyk00} (see also \cite{AKMW}).

It would be nice that the {\it weighted} blowups and blowdowns
 resulted from applying VGIT can be improved  to be just
{\it ordinary} blowups and blowdowns so that the un-weighted
factorization theorem may also have an easy and short GIT proof.
The problem boils down to the problem of what I call ``resolving
singular ${\mathbb C}^*$-action''. Namely, suppose we have a
$\C^*$-action over the union of two nonsingular open varieties $U
\cup V$ such that $\C^*$ acts freely with two smooth projective
quotients $U/\C^*$ and $V/\C^*$. Then we would like to find a
smooth equivariant compactification $W$ of $U \cup V$ such that
every isotropy subgroup of ${\mathbb C}^*$ on $W$ is either the
identity subgroup or the full group ${\mathbb C}^*$. Such an
action is called quasi-free. We do not know whether this can be
done.

 If this can be answered affirmatively,
 the method of {\bf \cite{HuKeel}} combined with \cite{DH} will imply that
 $f: X ---> X'$ can factorized as a sequence of
 (ordinary) blowup and blow downs along smooth centers. The key is
 the so-called Bialynicky-Birula decomposition theorem.


\subsection{VGIT and Weak Factorization Theorem for Projective Orbifolds}

A more exciting problem about VGIT $\Leftrightarrow$ Mori is to
generalize it to the category of Projective Varieties with Finite
Quotient Singularities. This is more compelling because, apart
from the intrinsic reasons, there is recently a surge of research
interests in orbifolds/stacks, for example, see the stringy
geometry and topology of orbifolds (\cite{RT}), among others.
These varieties occur naturally in many fields such as the Mirror
Symmetry Conjectures.

\begin{prob}
\label{stackyWF} Let $M ---> M'$ be a birational map between two
{\it projective} varieties with {\it finite quotient
singularities}. Decompose it as a composition of  explicit {\it
simple} birational transformations.
\end{prob}

To this end, we proved

\begin{thm} {\rm (--, {\bf \cite{Hu2004}})}
\label{thm:GIT}
 Let $\phi: X \rightarrow Y$ be a birational morphism between two  projective
varieties with at worst finite quotient singularities. Then there
is a smooth polarized projective $(\GL_n \times \C^*)$-variety
$(M, \cL)$ such that
\begin{enumerate}
\item $\cL$ is a very ample line bundle and admits two (general)
linearizations $\cL_1$ and $\cL_2$ with $M^{ss}(\cL_1) = M^s
(\cL_1)$ and $M^{ss}(\cL_2) = M^s (\cL_2)$. \item The geometric
quotient $M^s (\cL_1)/(\GL_n \times \C^*)$ is isomorphic to $X$
and the geometric quotient $M^s (\cL_2)/(\GL_n \times \C^*)$ is
isomorphic to $Y$. \item The two linearizations $\cL_1$ and
$\cL_2$ differ only by characters of the $\C^*$-factor, and
$\cL_1$ and $\cL_2$ underly the same linearization of the
$\GL_n$-factor. Let $\underline{\cL}$ be this underlying $\GL_n$-
linearization. Then we have $M^{ss} (\underline{\cL}) = M^s
(\underline{\cL})$.
\end{enumerate}
\end{thm}

The idea of the proof is as follows. First, we use two line
bundles over $X$, resulting from the pullback of an ample line
bundle over $Y$, to construct a $\C^*$-variety $Z$ (possibly very
singular) such that $X$ and $Y$ are identified with two geometric
quotients $Z^s(L_1)/\C^*$ and $Z^s(L_2)/\C^*$, respectively. Let
$V =Z^s(L_1) \cup Z^s(L_2)$. Then it has at most finite quotient
singularities because $X$ and $Y$ are such varieties and the
$\C^*$ action on $V$ is free. This implies that there is a
quasi-projective scheme $U$ with $\GL_n \times \C^*$ action such
that $V = U /\GL_n$ and $X$ and $Y$ are two geometric quotients of
$U$ by $\GL_n \times \C^*$. Then take an equivariant
compactification of $U$ and resolve singularities eqivariantly, we
finally arrive at our theorem above.

Next, we can verify that the action in the above theorem satisfies
the condition of the wall-crossing theorem (Theorem 4.2.7) of
\cite{DH}, hence as a consequence, we obtain

\begin{thm} {\rm (--,  {\bf \cite{Hu2004}})}
\label{thm:WFT}
 Let $X$ and $Y$ be two birational projective
varieties with at worst finite quotient singularities. Then $Y$
can be obtained from $X$ by a sequence of VGIT wall-crossing
birational transformations.
\end{thm}

In Theorem 1.2 of \cite{Hu2004}, as in this theorem, by a GIT
weighted blowup and a weighted blowdown, we mean the VGIT
wall-crossing birational transformation (or VGIT wall-crossing
flip for short) as described in Theorem 4.2.7 of \cite{DH}.

As a consequence of the proof, we also obtain the following
corollary to render the above suitable and useful for some
applications.

\begin{cor} {\rm (--, {\bf \cite{Hu2004}})}
\label{1torus}
 Let $X \rightarrow Y$ be a birational morphism between two
projective varieties with finite quotient singularities. Then $X$
and $Y$ can be realized as two geometric quotients of a projective
variety with finite quotient singularities by a ${\mathbb  C}^*$
action.
\end{cor}

The difference between Theorem \ref{thm:GIT} and Corollary
\ref{1torus} is that the former involves a reductive group with a
non-Abelian factor, while the latter uses only $\C^*$-action which
might be simpler from some perspectives, but for this, one has to
sacrifice the smoothness of the ambient variety.

Establishing all the above, a  natural next question is

\begin{prob}
\label{invariants}
Describe how orbifold cohomology, orbifold Chow group,
and other orbifold invariants change when crossing a wall.
\end{prob}

Obviously, a deeper problem is to compute how orbifold quantum
cohomology changes when crossing a codimensional one wall (see
\cite{Kai}, \cite{RT}, and methods of \cite{Li1}, \cite{Li2}),
however, even at the ``classical'' level (without quantum
corrections), this problem is already very interesting and highly
non-trivial.

To this end, it may be  useful to bring a new cohomology theory
into the picture: the orbifold intersection cohomology for
singular orbispces, for which the very singular GIT quotient
defined by a polarization on the wall, which is also the common
blowdown of two adjacent smooth GIT quotient orbifolds, might be a
natural example. This is, of course, interesting in its own right.
The total orbifold intersection cohomology may be defined as the
total intersection cohomology of the twisted sectors, or
stack-theoretically, ``inertia stack''. To make it truly useful
for our purpose, a BBD type decomposition theorem for orbifold
cohomology may be needed (\cite{BBD}). Also to better track how
orbifold invariants change, we need

\begin{prob}
\label{inertia} Describe explicitly the twisted sectors of  GIT
quotient orbifolds  using the raw data from the group action, and
describe explicitly
how it changes when crossing a wall.
\end{prob}

The last two problems are, in fact, already interesting in the
settings of \cite{HuMIT},  in that even if the ambient variety is
smooth, the GIT quotients may still have quotient singularities,
and then, according to stringy geometry, orbifold invariants are
the right ones we should study. In other words, we expect to
``refine'', in the case of the presence of (finite) quotient
singularities, some of our intersection cohomology calculations in
\cite{HuMIT},  in terms of the new orbifold invariants.


\section{Chow Quotients: Perturbation-Translation-Specialization}

GIT quotient (or stability) can be very difficulty to handle in
practice. But Chow quotient in general is even harder. The
difficulty of Chow quotient is that it is in general very
singular. However, understanding Chow quotient well may lead to
some significant applications (\cite{Hu2005}).

\subsection{Definition}

Consider a reductive algebraic group action on a projective
variety over the field of complex numbers
$$G \times X \to X.$$
The Chow quotient of this action is defined as follows.

There is a small open subset, $U \subset X$, such that
$\overline{G \cdot x}$ represents the same Chow cycle $\delta$ for
all $x \in U$. Let Chow$_\delta(X)$ be the component of the Chow
variety of $X$ containing $\delta$. Then there is an embedding
$$\iota: U/G \to \hbox{Chow}_\delta (X)$$
$$ [G \cdot x] \to [\overline{G \cdot x}] \in \hbox{Chow}_\delta (X).$$
The Chow quotient, denoted by $X/\!/^{ch}G$, is defined to be 
the closure of $\iota (U/G)$. This definition is independent of
the choice of the open subset $U$. That is, the Chow quotient is
canonical.

\subsection{Chow Family}

There is a "universal" family over the Chow quotient. Let $$F
\subset X \times (X/\!/^{ch}G)
$$ be the family of algebraic cycles over the Chow quotient
$X/\!/^{ch}G$ defined by
$$F= \{ (x, Z) \in X \times (X/\!/^{ch}G) | x \in Z\}.$$
Then, we have a diagram
\begin{equation*}
\begin{CD}
F @>{\ev}>> X\\
@V{f}VV  \\
X/\!/^{ch}G
\end{CD}
\end{equation*}
where $\ev$ and $f$ are  the projections to the first and second
factor, respectively.  For any point $q \in X/\!/^{ch}G $, we will
call the fiber, $f^{-1}(q)$, the Chow fiber over the point $q$.
Sometimes we identify $f^{-1}(q)$ with its embedding image
$\ev(f^{-1}(q))$ in $X$.

\subsection{Perturbation-Translation-Specialization}

As we mentioned before, the Chow quotient approach to moduli
problems is different and may be harder than GIT, but it has some
advantages and potentially has several significant applications.
Understanding the boundary of Chow quotient is the key to this
approach.

 In our paper \cite{Hu2003b},
we showed that for a torus $G$ action on a smooth projective
variety $X$, the boundaries of the Chow quotients, or
equivalently, the special fibers of the Chow family,  admit a
computable characterization. This
 is very   desirable for applications. Roughly, we prove that

\begin{thm} {\rm (--, {\bf \cite{Hu2003b}})}
\label{chowfiberAG}
 Two points $x$ and $y$ of $X$, with
$$\dim G \cdot x = \dim G \cdot y = \dim G,$$
are in the same Chow fiber if and only if $x$ can be perturbed (to
general positions), translated along $G$-orbits (to positions
close to $y$), and then specialized to the point $y$. (See Figure
3 for an illustration.)
\end{thm}

\bigskip

\begin{picture}(12, 11)
\put(3,1){ \psfig{figure=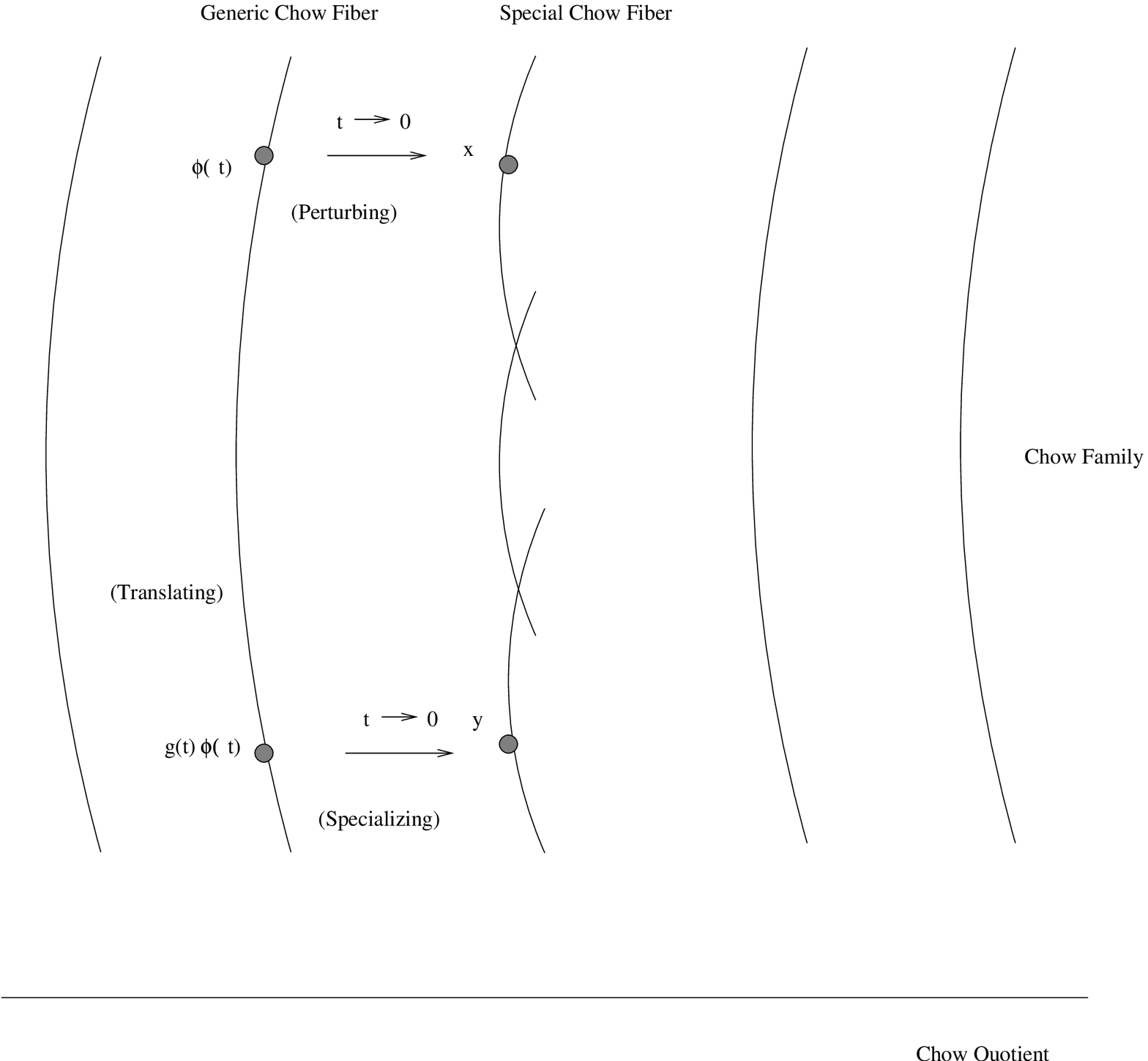,height = 9.5cm,width = 10.5cm} }
\end{picture}

\centerline{Figure 3.  Perturbing-Translating-Specializing}




This relation, which we call {\it
perturbing-translating-specializing} (P.T.S.) relation,  was first
discovered and studied intensively by Neretin in the case of
symmetric spaces. In addition to the above, we also proved in
\cite{Hu2003b} that over the field of complex numbers, the Chow
quotient admits
  symplectic and other topological interpretations, namely,
\begin{itemize}
\item symplectically,
  the moduli spaces of stable orbits with prescribed momentum charges; and
  \item topologically,
  the moduli space of stable action-manifolds.
\end{itemize}

Moduli spaces of what I call stable polygons are symplectic
interpretations of the Chow quotient $\overline{M}_{0,n}$, the
moduli space of stable $n$-pointed rational curves. For detailed
account of this and beyond, see \cite{Hu99}.

An upshot of PTS relation is that, comparing to the nondescriptive
definition of special Chow fibers, it is computable, and thus
provides some much needed information on boundary cycles of the
Chow quotient. As an application, we have appllied Theorem
\ref{chowfiberAG} to the case of point configurations on ${\Bbb
P}^n$ ($n >1$), and propose a {\it geometric interpretation} of
the Chow quotients of $({\Bbb P}^n)^m$ (equivalently, the Chow
quotients of higher Grassmannians).

In \cite{Lafforgue}, Lafforgue provides some toric
compactifications of matroid strata of Grassmannians. His
compactifications are reducible in general. Our P.T.S theorem may
be applied to single out the main irreducible components.

The PTS relation can be used to explain more naturally the moduli
space  $\overline{M}_{0,n}$ as the Chow quotient of $({\mathbb
P}^1)^n$ (\cite{Hu2004b}). It may also be useful for its toric
degeneration (\cite{Hu2005a}).

\section{Birational Geometry of HyperK\"ahler Varieties}

HyperK\"ahler Varieties form a sepcial and very important class of
Calabi-Yau manifolds. Its birational geometry has received some
considerable attentions recently (\cite{BHL}, \cite{HY},
\cite{Fu}, \cite{WW},among others). Being very restrictive in
structures, birtaional transformations among HyperK\"ahler
Varieties are among simplest kinds. We will be mainly interested
in projective symplectic varieties. This means a nonsingular
projective variety $X$ of dimension $2n$, together with a
holomorphic two form $\omega$ such that $\omega^n \ne 0$.

Let $f: X ---> Y$ be a birational map between two symplectic
varieties. By the weak Factorization Theorem, it can be factorized
as blowup and blowdowns. However, we have little control on the
varieties in between, in particular, we do not know whether they
are symplectic or not. This causes difficulty to keep track the
changes of symplectic varieties. Hence, special efforts are need
to factorize $f$ in the category of symplectic varieties.

\subsection{Factorizations in Dimension 4}

The simplest known  birational transformation among symplectic
varieties are the so-called Mukai elementary transformations
(MET). For a symplectic 4-fold $X$, this means the following. Let
${\mathbb P}^2$ be embedded in $X$, then it is a Lagrangian
subvariety, hence its normal bundle is isomorphic to its cotangent
bundle. Blowup $X$ along this ${\mathbb P}^2$, then its
exceptional divisor $D$ is the incident variety
$$D = \{(p,l) \in {\mathbb P}^2 \times \check{{\mathbb P}}^2 | \; p \in l\}$$
where  $\check{{\mathbb P}}^2$ is the space of lines in ${\mathbb
P}^2$. Mukai showed that  it can be contracted along another
ruling, producing another symplectic variety $X'$. This process is
called a MET.

The same can be done in general, when $X$ contains a projective
bundle whose rank is the same as its codimension.

In dimension 4, we have the following structure theorem.

\begin{thm}
{\rm (Burns-Luo-Hu, {\bf \cite{BHL}})} Let $f: X ---> Y$ be a
birational map between two symplectic varieties of dimension 4.
Assume that every component of the degenerate locus of $f$ is
normal. Then $f: X
---> Y$ can be factorized as a sequence of Mukai elementary
transformations.
\end{thm}

The normality condition in the theorem was later removed by
Wierzba and Wisniewski (\cite{WW}).

\subsection{Symplectic Contractions in Higher Dimensions}

The general case is much more difficult. As an important initial
step, Yau and I obtained the following structure theorems on small
symplectic contractions.

\begin{thm} {\rm (Hu-Yau, {\bf \cite{HY}})}
Let $X$ be a smooth projective symplectic variety and $\pi: X \to
Z$ be a small contraction. Let $B$ be an irreducible component of
the degenerate locus of $\pi$ and $F$ is a generic fiber of the
restricted map $\pi: B \rightarrow S= \pi(B)$. Then we have
\begin{enumerate}
\item $T_b F = (T_b B)^\perp$ for a generic smooth point $b \in F$
and the orthogonal space is taken with respect to the symplectic
form $\omega$. \item the inclusion $j: B \hookrightarrow X$ is a
coisotropic embedding. \item the null foliation of $\omega|_B$
coincides the generic fibers of $\pi: B \rightarrow S$.
\end{enumerate}
\end{thm}

As a special case, we showed

\begin{thm} {\rm (Hu-Yau, {\bf \cite{HY}})}
Assume that $\pi: B \rightarrow S$ is a smooth fiberation, then it
must be a projective bundle with rank $r = \codim B$. In
particular, the Mukai transformation can be performed along this
projective bundle to get another symplectic variety.
\end{thm}

\begin{cor} {\rm (Hu-Yau, {\bf \cite{HY}})}
If $B$ is smooth and can be contracted to a point, then $B \cong
{\mathbb P}^n$, in particular, it is a Lagrangian. In fact,
whenever $B$ can be contracted to a point, smooth or not, it must
be a Lagrangian.
\end{cor}

\begin{prop} {\rm (Hu-Yau, {\bf \cite{HY}})}
Any small contraction $\pi: X \to Z$ must be perverse semi-small
in the sense that there is a stratification $Z= \cup Z_\alpha$
such that $$\dim \pi^{-1}(z) \le {1 \over 2} \codim Z_\alpha$$ for
all $z \in Z_\alpha$ and non-open strata $Z_\alpha$.
\end{prop}

When all the inequalities are equalities, the map is said to be
strictly semi-small. We believe

\begin{conj} {\rm (Hu-Yau, {\bf \cite{HY}})} Any small contraction $\pi: X \to Z$
must be perverse strictly semi-small.
\end{conj}

This helps to convince us the following general conjecture

\begin{conj} {\rm (Hu-Yau, {\bf \cite{HY}})}
Let $f: X ---> Y$ be a birational map between two symplectic
varieties. Then after removing subvarieties of codimension greater
than 2, $X$ and $Y$ are related by a sequence of Mukai elementary
transformations.
\end{conj}

This conjecture was recently confirmed by B. Fu (\cite{Fu}) in the
case of the closures of nilpotent orbits.

We remark that some of our results in this section overlap with
the work of \cite{miyaoka}

\bigskip

{\sl Acknowledgements.} This article is by no means intended to be
a complete account of all the related development. Therefore I
apologize for the omission of any relevant work and references. I
would like to thank ICCM-2004 committee for the invitation and IMS
of the Chinese University of Hong Kong for financial support.
 Also, thanks to Dan Abramovich, J\'anos Koll\'ar and Hsian-hua Tseng for
useful mathematical comments. Special thanks go to Professor S.-T.
Yau for his valuable advice and support. Some of the research in
this article are partially supported by NSF and NSA.



\begin{thebibliography}{10}

\bibitem{AB} M. Atiyah and R. Bott,
{\em The moment map and equivariant cohomology.} Topology  23
(1984), no. 1, 1--28.  MR0721448.
\bibitem{AKMW} Dan Abramovich, Kalle Karu, Kenji Matsuki, Jaroslaw Wlodarczyk
 {\it  Torification and Factorization of Birational Maps},
 J. Amer. Math. Soc.  15  (2002),  no. 3, 531--572. MR1896232, Zbl1032.14003.
\bibitem{BassH} H. Bass and W. Haboush,
{\it Linearizing certain reductive group actions,}
Trans. Amer. Math. Soc. {\bf 292} (1985), no. 2, 463--482.
\bibitem{BP} M Brion and C. Procesi,
{\em Action d'un tore dans une vari\'et\'e,} in "Operator algebra,
unitary representations, enveloping algebras, and invariant
theory", Progress in Mathematics 192 (1990), Birkh\"auser,
509-539.
\bibitem{Kai} K. Behrend,
{\em Cohomology of Stacks,} University of British Colombia (2002).
\bibitem{BBD} A. Beilinson, J. Bernstein, and P. Deligne,
{\em Faisceaux Pervers.}  In: Analyse et topologie sur les espaces
singuliers, I. Ast\'erisque \#100, Soc. Math. de France (1983).
\bibitem{BH} H. Boden and Y. Hu,
{\em Variation of Moduli of Parabolic Bundles,} Mathematische
Annalen Vol. 301 No. 3 (1995), 539 -- 559.
\bibitem{BHL} D. Burns, Y. Hu and T. Luo,
{\em HyperK\"ahler manifolds and birational transformations in
dimension 4,} In: Vector bundles and representation theory,
Contemporary Mathematics  Volume {\bf 322}  (2003), 141 -- 149.
\bibitem{miyaoka} K. Cho, Y. Miyaoka, and N. Sheperd-Barron,
{\em Characterizations of projective space and applications to
complex symplectic manifolds.}  Higher dimensional birational
geometry (Kyoto, 1997),  1--88, Adv. Stud. Pure Math., 35, Math.
Soc. Japan, Tokyo, 2002. MR1929792
\bibitem{Cox} D. Cox,
{\em The homogeneous coordinate ring of a toric variety. }
 J. Algebraic Geom.  {\bf 4}  (1995),  no. 1, 17--50.
\bibitem{cutkosky} S. D. Cutkosky,
{\em  Strong Toroidalization of Birational Morphisms of 3-Folds.}
math.AG/0412497
\bibitem{DH} I. Dolgachev and Y. Hu,
{\em Variation of Geometric Invariant Theory Quotients,} Publ.
Math. I.H.E.S. {\bf 87} (1998), 5 -- 56. MR1659282, Zbl
1001.14018.
\bibitem{Donaldson} S. Donaldson,
{\em Scalar curvature and projective embeddings. I,}
 J. Differential  Geom. {\bf 59} (2001), no. 3, 479--522.
\bibitem{Donaldson2002} S. Donaldson,
{\em Scalar curvature and stability of toric varieties .}
J. Differential  Geom. Volume {\bf 62} (2002), 289--349.
\bibitem{EHKV} D. Edidin, B. Hassett, A. Kresch, A.
Vistoli, {\em Brauer Groups and Quotient Stacks}.
 Amer. J. Math.  123  (2001),  no. 4, 761--777. MR1844577, Zbl 1036.14001.
 \bibitem{GS} V. Guillemin and S. Sternberg,
 {\em Birational equivalence in the symplectic category. } Invent. Math.  97  (1989),  no. 3, 485--522.
 MR1005004.
\bibitem{GM} M. Goresky and R. MacPherson,
{\em On the topology of torus actions,} in ``Algebraic Groups''
(1985).
\bibitem{Faltings} G. Faltings and G. W\"ustholz,
{\em Diophantine approximations on projective spaces,} Invent.
Math {\bf 116} (1994), 109--138.
\bibitem{FH} P. Foth and Y. Hu,
{\em Toric degenerations of weight varieties and
applicationsmath,} Travaux Mathématiques (2005), math.AG/0406329.
\bibitem{Fu} B. Fu,
{\em Birational geometry in codimension 2 of symplectic
resolutions.} math.AG/0409224
\bibitem{Giv} A. Givental,
{\em Equivariant Gromov-Witten invariants.} Internat. Math. Res.
Notices  1996,  no. 13, 613--663. MR1408320
\bibitem{Ha} P. Hacking,
{\em Compact moduli of hyperplane arrangements,}
    math.AG/0310479.
\bibitem{HuMITthesis} Y. Hu,
{\em Geometry and topology of quotient varieties,} Ph.d Thesis,
MIT , 1991.
\bibitem{HuMIT} Y. Hu,
{\em Geometry and topology of quotient varieties of torus
actions,} Duke Math. Journal {\bf 68} (1992), 151 -- 183.
\bibitem{Hu96} Y. Hu,
{\em Relative geometric invariant theory and universal moduli
spaces.}  Internat. J. Math. {\bf 7 }(1996), no. 2, 151--181.
MR1382720, Zbl 0889.14005.
\bibitem{Hu99} Y. Hu,
 {\em Moduli Spaces of Stable Polygons and Symplectic Structures on
$\overline{M}_{0,n}$,} Compositio Mathematica Vol. 118 (1999), 159
-- 187.
\bibitem{Hu2003a} Y. Hu,
{\em  Stable Configurations of Linear Subspaces and Vector
Bundles.}  Quarterly Journal of Pure and Applied Mathematics
Volume 1 No. 1 (2005), 127 -- 164.
\bibitem{Hu2003b} Y. Hu,
{\em  Topological Aspects of Chow Quotients.} Preprint, submitted
(2003).  math.AG/0308027
\bibitem{Hu2004} Y. Hu,
{\em Factorization Theorem for Projective Varieties with Finite
Quotient Singularities,}  Journal of Differential Geometry Volume
{\bf 68} (2004), 587 -- 593.
\bibitem{Hu2004b} Y. Hu,
{\em $\overline{M}_{0,n}$ as Chow Quotient Revisited,} In
preparation.
\bibitem{Hu2005a} Y. Hu,
{\em Toric Degeneration of $\overline{M}_{0,n}$.} In preparation.
\bibitem{Hu2005} Y. Hu,
{\em } In preparation.
\bibitem{HuKeel2} Y. Hu and S. Keel,
{\em Mori dream spaces and GIT,}
Dedicated to William Fulton on the occasion of his 60th birthday.
 Michigan Math. J.  {\bf 48}  (2000), 331--348. MR1786494, Zbl pre01700876.
\bibitem{HuKeel} Y. Hu and S. Keel,
{\em A GIT proof of the weighted weak factorization theorem.}
unpublished short notes, available at Math ArXiv, math.AG/9904146.
\bibitem{huyb2} D. Huybrechts,
{\em Compact HyperK\"ahler  manifolds: basic results,} Invent.
Math.  {\bf 135} (1999), 63-113.
\bibitem{HL} Y. Hu and W. Li,
{\em Variation of the Gieseker and Uhlenbeck Compactifications,}
International Journal of Mathematics  Vol. 6 No. 3 (1995), 397 --
418.
\bibitem{HY} Y. Hu and S.-T. Yau,
{\em HyperK\"hler Manifolds and Birational Transformations,}
 Adv. Theor. Math. Phys. Volume {\bf 6} Number 3 (2002), 557-576.
 \bibitem{HLY} Y. Hu, C.-H. Liu, and S.-T. Yau,
 {\em Toric Morphisms and Fibrations of Toric Calabi-Yau
hypersufaces,} Adv. Theor. Math. Phys. Volume 6 Number 3 (May
2002), 457 -- 506.
\bibitem{Ka} M. Kapranov,
{\em Chow quotients of Grassmannian,} I.M. Gelfand Seminar
Collection, 29--110, Adv. Soviet Math., 16, Part 2, Amer. Math.
Soc., Providence, RI, 1993.
\bibitem{KSZ} M. Kapranov, B. Sturmfels and A. Zelevinsky,
{\it Quotients of toric varieties}, Math. Ann. 290 (1991), no. 4,
643--655.
\bibitem{KM} S. Keel and S. Mori,
{\em  Quotients by groupoids.}  Ann. of Math. (2)  145  (1997),
no. 1, 193--213.
\bibitem{KT} S. Keel and J. Tevelev,
{\em  Chow Quotients of Grassmannians II,} math.AG/0401159.
\bibitem{KNess} G. Kempf and L. Ness,
{\em The length of vectors in representation spaces. } Algebraic
geometry (Proc. Summer Meeting, Univ. Copenhagen, Copenhagen,
1978), pp. 233--243, Lecture Notes in Math., 732, Springer,
Berlin, 1979. MR0555701.
\bibitem{Kirwan84} F. Kirwan,
{\em Cohomology of quotients in symplectic and algebraic
geometry.} Mathematical Notes, 31. Princeton University Press,
Princeton, NJ, 1984.  MR0766741
\bibitem{Kirwan84} F. Kirwan,
{\em Partial desingularisations of quotients of nonsingular
varieties and their Betti numbers.}  Ann. of Math. (2)  122
(1985), no. 1, 41--85.  MR0799252.
\bibitem{Kly} A. Klyachko,
{\em Stable bundles, representation theory and Hermitian operators,}
Selecta Math. {\bf 4} (1998), 419-445.
\bibitem{Ko97} J. Koll\'ar,
{\em Quotient spaces modulo algebraic groups.}  Ann. of Math. (2)
145  (1997),  no. 1, 33--79.
\bibitem{km} J. Koll\'ar and S. Mori,
{\em Birational geometry of algebraic varieties,}
Cambridge University Press, 1998.
\bibitem{Lafforgue} L. Lafforgue,
{\em Chirurgie des grassmanniennes.}  (French) [Surgery on
Grassmannians] CRM Monograph Series, 19. American Mathematical
Society, Providence, RI, 2003.
\bibitem{Lafforgue99} L. Lafforgue,
 {\em Pavages des simplexes, schémas de graphes recollés et compactification
 des ${\rm PGL}\sp {n+1}\sb r/{\rm PGL}\sb r$}, Invent. Math.  136  (1999),  no. 1, 233--271.
 Erratum:   Invent. Math.  145  (2001),  no. 3, 619--620.
\bibitem{Li1} Jun Li, {\em Stable morphisms to singular schemes
and relative stable morphisms}, J. Diff.\ Geom.\ {\bf 57} (2001),
no.\ 3, 509--578.
\bibitem{Li2} Jun Li, {\em A degeneration formula of GW-invariants},
J.\ Diff.\ Geom.\ {\bf 60} (2002), no.\ 2, 199--293.
\bibitem{Yau1} B, Liang, K. Liu and S. T. Yau,
{\em Mirror Principles, I} Asian J. of Math. {\bf 1} (1997).   no.
4, 729--763. MR1621573
\bibitem{Yau2} B, Lian, K. Liu and S. T. Yau,
{\em Mirror Principles, II,  III,} Asian J. of Math. (1999) no. 1,
109--146 and   no. 4, 771--800,   MR1701925 and MR1797578
 \bibitem{Yau3} B, Lian, K. Liu and S. T. Yau,
 {\em Mirror principle. IV }  Surveys in differential geometry,
 433--474, Surv. Differ. Geom., VII, Int. Press, Somerville, MA, 2000.
\bibitem{Luo} H.-Z. Luo,
{\em Geometric criterion for Gieseker-Mumford stability of polarized manifolds,}
Journal of Diff. Geom. {\bf 49} (1998), no. 3, 577--599.
\bibitem{GIT} D. Mumford,
{\em Geometric Invariant Theory,} 1962. Berlin, New York.
MR1304906, Zbl 0797.14004.
\bibitem{Ness} L. Ness,
{\em Stratification of the null cone via the moment map.} With an
appendix by David Mumford. Amer. J. Math. 106 (1984), no. 6,
1281--1329. MR0765581
\bibitem{Pandharipande94a} R. Pandharipande,
{\em A compactification over $\overline{M_g}$ of the universal
moduli space of slope-semistable vector bundles,}  J. Amer. Math.
Soc.  9  (1996),  no. 2, 425--471.  MR1308406
\bibitem{Reichstein89} Z. Reichstein,
{\em  Stability and equivariant maps},
 Inven. Math. {\bf 96} (1989),
 349-383.
\bibitem{PhS} D.H. Phong and J. Sturm,
{\em Stability, Energy Functions, and K\"ahler-Einstein Metrics,}
arXiv:math.DG/0203254.
\bibitem{RT1} J. Ross and R. P. Thomas,
{\em A study of the Hilbert-Mumford criterion for the stability of
projective varieties. } math.AG/0412519
\bibitem{RT2} J. Ross and R. P. Thomas,
{\em An obstruction to the existence of constant scalar curvature
K\"ahler metrics}, math.DG/0412518
\bibitem{RT} Y. Ruan,
{\em  Stringy geometry and topology of orbifolds. } Symposium in
Honor of C. H. Clemens (Salt Lake City, UT, 2000), 187--233,
Contemp. Math., 312, Amer. Math. Soc., Providence, RI, 2002.
\bibitem{Simpson94} C. Simpson,
{\em Moduli spaces of representations of the fundamental group of
a smooth projective variety, I and II},  Inst. Hautes \'Etudes
Sci. Publ. Math.  No. 79 (1994), 47--129.  MR1307297.  Inst.
Hautes \'Etudes Sci. Publ. Math.  No. 80 (1994), 5--79 (1995).
MR1320603.
\bibitem{Th} M. Thaddeus,
{\em Geometric invariant theory and flips,}  J. Amer. Math. Soc. 9
(1996),  no. 3, 691--723. MR1333296, Zbl 0874.14042.
\bibitem{tian} G. Tian,
{\em K\"ahler-Einstein metrics on algebraic manifolds,}
LNM {\bf 1646}, Spinger, Berlin, 1996.
\bibitem{Tian97} G. Tian,
{\em K\"ahler-Einstein metrics with positive scalar curvature.}
Invent. Math. {\bf 130} (1997), 1--39.
\bibitem{Totaro} B. Totaro,
{\em Tensor products of semistables are semistable,}
Geometry and Analysis on complex manifolds, World Sci. Publ. 1994, pp. 242-250.
\bibitem{Wang} X. W. Wang,
{\em Balanced point, stability and vector bundles over projective
manifolds.}
 Math. Res. Lett.  {\bf 9}  (2002),  no. 2-3, 393--411.
\bibitem{Viehweg95} E. Vieweg,
{\em Quasi-projective moduli of polarized manifolds,}
Springer-Verlag, 1995.
 \bibitem{WW} J. Wierzba and J. Wi\'sniewski,
 {\em Small contractions of symplectic 4-folds. }
 Duke Math. J.  120  (2003),  no. 1, 65--95. MR2010734.
 \bibitem{Wlodarczyk00}
J. Wlodarczyk, {\em Birational cobordisms and factorization of
birational maps.}  J. Algebraic Geom.  9  (2000), no. 3, 425--449.
MR1752010, Zbl 1010.14002.
\bibitem{Wlodarczyk03}
J. Wlodarczyk, {\em Toroidal varieties and the weak factorization
theorem.} Invent. Math.  154  (2003),  no. 2, 223--331. MR2013783,
Zbl pre02021105.
\bibitem{Wlodarczyk97}
J. Wlodarczyk, {\em Decomposition of birational toric maps in
blow-ups \& blow-downs.}  Trans. Amer. Math. Soc. 349 (1997), no.
1, 373--411. MR1370654
\bibitem{Yau} S.-T. Yau,
{\em Open problems in geometry,} Proc. Symp. Pure Math. {\bf 54}
(1993) 1-28.
\bibitem{Yau87} S.-T. Yau,
{\it Mathematical aspects of string theory.} Proceedings of the
conference held at the University of California, San Diego,
California, July 21--August 1, 1986. Edited by S.-T. Yau. Advanced
Series in Mathematical Physics, 1. World Scientific Publishing
Co., Singapore, 1987. MR0915812
\bibitem{zhang} S. Zhang,
{\em Heights and reductions of semi-stable varieties,}
Compositio Mathematica {\bf 104} (1996) 77-105.

\end{thebibliography}
\end{document}